\documentclass[11pt,tbtags,reqno]{article}
\usepackage[margin=0.90in]{geometry}
\usepackage{amssymb}
\usepackage{amsmath}
\usepackage{amsthm}
\usepackage[swedish, english]{babel}
\usepackage{psfrag}   
\usepackage{graphicx}
\usepackage[retainorgcmds]{IEEEtrantools}
\usepackage{bbm}
\usepackage[square, comma, numbers, sort&compress]{natbib}
\usepackage{enumerate}
\usepackage{soul}
\usepackage{hyperref}
\usepackage[colorinlistoftodos, textsize=tiny, backgroundcolor=yellow
, linecolor=green]{todonotes}
\usepackage{sidecap}
\usepackage{float}
\usepackage{leftidx}

\usepackage{bm}

\newcommand{\E}{\mathbb E}
\newcommand{\R}{\mathbb R}
\newcommand{\N}{\mathbb N}

\newcommand{\ep}{\varepsilon}

\def\P{{\mathbb P}}

\numberwithin{equation}{section}
\theoremstyle{plain}

\newtheorem{theorem}{Theorem}[section]
\newtheorem{lemma}[theorem]{Lemma}

\newtheorem{proposition}[theorem]{Proposition}

\theoremstyle{definition}
\theoremstyle{example}

\begin{document}

\title{\textsc{
A Strong Law of Large Numbers for Positive Random Variables} \thanks{~ We are deeply grateful to J\'anos \textsc{Koml\'os}, who went over the entire manuscript with the magnifying glass and offered line-by-line criticism and wisdom. We thank   Daniel \textsc{Ocone},  Albert \textsc{Shiryaev} for invaluable advice; and Richard \textsc{Groenewald}, Tomoyuki \textsc{Ichiba}, Kostas \textsc{Kardaras},  Tze-Leung \textsc{Lai},  Kasper \textsc{Larsen}, Ayeong \textsc{Lee}, Emily \textsc{Sergel}, Nathan \textsc{Soedjak}   for careful readings     and   suggestions.}
}

\author{  
\textsc{Ioannis Karatzas} \thanks{~
Department of Mathematics,  Columbia University, New York, NY 10027 (e-mail: {\it ik1@columbia.edu}). Support from the National Science Foundation under Grant  DMS-20-04977 is gratefully acknowledged.  
}  
 \and
\textsc{Walter Schachermayer}                \thanks{~ 
Faculty of Mathematics, University of Vienna, Oskar-Morgenstern-Platz 1, 1090 Vienna, Austria (email: {\it walter.schachermayer@univie.ac.at}). Support from the  Austrian Science Fund (FWF) under
 grant P-28861,  and by the Vienna Science and  Technology Fund (WWTF) through project MA16-021, is gratefully acknowledged.   
          }
                                      }

\maketitle
\begin{abstract}
\noindent
In the spirit of the famous \textsc{Koml\'os} (1967) theorem, every sequence  of nonnegative, measurable functions $\{ f_n \}_{n \in \N}$    on a probability space, contains a subsequence which---along with all its subsequences---converges a.e.\,\,in \textsc{Ces\`aro} mean to some   measurable    $f_* : \Omega \to [0, \infty]$.   This result of \textsc{von\,Weizs\"acker} (2004) is proved here using a new methodology  and elementary tools; these  sharpen  also a  theorem of \textsc{Delbaen \& Schachermayer} (1994),   replacing general convex combinations by \textsc{Ces\`aro}  means. 
\end{abstract}

\noindent
 {\sl AMS  2020 Subject Classification:}  Primary 60A10, 60F15; Secondary 60G42, 60G46.

\noindent
 {\sl Keywords:}  Strong law of large numbers, hereditary convergence,   partition of unity


\section{Introduction}
\label{sec1}


 On a probability space $(\Omega, \mathcal{F}, \mathbb{P})$,   consider    real-valued   measurable functions $f_1, f_2, \cdots \,.$ If these are independent and have  the same distribution with     $\E ( | f_1|) < \infty\,$, the celebrated \textsc{Kolmogorov}  strong law of large numbers (\cite{K2};\,\cite{KAN};\,\cite{Du}, p.\,73) states that the ``sample average" $ \, ( f_1 + \cdots + f_N) / N\,$ converges $\mathbb{P}-$a.e.\,to the ``ensemble average" $\,\E (f_1) = \int_\Omega f_1 \, \mathrm{d} \mathbb{P}\,,$ as $ N \to \infty$.  More generally,  if $f_n (\omega) = f \big( T^{n-1}(\omega) \big), \, n \ge 2, \, \omega \in \Omega$ are the images of an integrable function $f_1: \Omega \to \R$ along the orbit  of successive actions of a measure-preserving transformation $T:   \Omega \to \Omega\,,$ then the  above sample average converges $\mathbb{P}-$a.e.\,to the conditional expectation  $f_* = \E ( f_1 |{\cal I})$  of $f_1$ given the $\sigma-$algebra   ${\cal I}$  of   $T-$invariant sets,  by the \textsc{Birkhoff} pointwise ergodic theorem (\cite{Du}, p.\,333).

A deep  result of \textsc{Koml\'os} \cite{K}, already   55 years  old but always very striking, says that such ``stabilization via  averaging" occurs within {\it any} sequence $f_1, f_2, \cdots \,$ of measurable, real-valued  functions with  $\, \sup_{n \in \N} \E ( | f_n|) < \infty\,.$ More precisely, there exist  then an integrable function $f_*$ and a subsequence $ \{ f_{n_k}  \}_{k \in \N} $ such that   $ \, ( f_{n_1} + \cdots + f_{n_K}) / K\,$ converges to $f_*\,$, $\,\mathbb{P}-$a.e.\,\,as $ K \to \infty$; and   the same is true for any further subsequence  of this $ \{ f_{n_k}  \}_{k \in \N}\,.$  

 \bigskip
 \newpage
 
This result inspired further path-breaking work in probability  theory  (\cite{G},\,\cite{Ch2},\,\cite{Ch3})  culminating with   \\    \textsc{Aldous}      
 (1977),  where exchangeability plays   a crucial r\^ole.    It, and its ramifications \cite{DS1},\,\cite{DS2} involving forward convex combinations, have been   very useful in   the field of convex optimization;   more generally, when one seeks objects with specific properties, and tries to ascertain
 their existence using weak compactness arguments. Stochastic control, optimal stopping and      hypothesis testing   are   examples of the former (e.g.,\,\cite{KS},\,\cite{KW},\,\cite{CK},\,\cite{KZ},\,\cite{LZ});    the   \textsc{Doob-Meyer}  and \textsc{Bichteler-Dellacherie}  theorems in stochastic analysis provide instances    of the latter (e.g.,\,\cite{J},\,\cite{BSV1},\,\cite{BSV2}).  
 
 We develop here a very simple   argument for the \textsc{Koml\'os} theorem, in the important special case of nonnegative   $f_1, f_2, \cdots \,$ treated by \textsc{von\,Weizs\"acker} (2004).  The argument    dispenses   with   boundedness in $ \mathbb{L}^1$,  at the cost of   allowing  the function $f_*$ to take  infinite values. 

\section{Background}
\label{sec2}

We place ourselves   on a given, fixed probability space $(\Omega, \mathcal{F}, \mathbb{P})$, and consider  a sequence $f_1, f_2, \cdots \,$ of measurable, real-valued functions  defined on it. 
We say that this sequence   {\it converges  hereditarily in \textsc{Ces\`aro} mean} to some measurable $f_*:\Omega \to \mathbb{R} \cup \{\pm \infty\}$, and write 
$ 
f_n \xrightarrow[n \to \infty]{hC}f_*\,,~~ \mathbb{P} -\hbox{a.e.,} 
$  
if, for  {\it every} subsequence $\big\{f_{n_k} \big\}_{k \in \mathbb{N}}$ of the original sequence, we have 
\begin{equation}
\label{1}
\lim_{K\to \infty} \frac{1}{K} \sum_{k=1}^K f_{n_k} = f_*\,,\qquad \mathbb{P} -\hbox{a.e.}
\end{equation}
 Clearly then,   every other such sequence $g_1, g_2, \cdots \,$ which is 
 {\it equivalent} to   $f_1, f_2, \cdots \,,$   in the sense of $\, \sum_{n \in \N} \mathbb{P}( f_n \neq g_n) < \infty\,$ (cf.\,\cite{KAN}), also has this property. 

\smallskip
In 1967,  \textsc{Koml\'os}   proved the following  remarkable result. The argument  in \cite{K} is very clear, but also  long and  quite   involved. Simpler proofs and extensions have appeared since (e.g.,\,\cite{S},\,\cite{T};\,\cite{B}). 

\begin{theorem} 
[\textsc{Koml\'os} (1967)]
\label{Kom}
If the sequence $\{f_n \}_{n \in \mathbb{N}}$ is bounded  in $ \mathbb{L}^1,$  i.e.,
  $\sup_{n \in \mathbb{N}} \mathbb{E} (|f_n|) < \infty\,$ holds,     there exist an integrable   $f_*:\Omega \to \mathbb{R}$  and a subsequence $\big\{f_{n_k}\big\}_{k \in \mathbb{N}}$ of $\{f_n \}_{n \in \mathbb{N}}\,,$ which converges hereditarily in \textsc{Ces\`aro} mean to $f_*\,:$   
\begin{equation}
\label{02}
f_{n_k} \xrightarrow[k \to \infty]{hC}f_*\,, \qquad \mathbb{P}-\hbox{a.e.} 
\end{equation}
\end{theorem}

This result  was motivated by  an earlier one, Theorem \ref{Rev} right below.  For the convenience of the reader, we provide in \S \,\ref{sec5f} a simple    proof  (in the manner of\,\cite{Ch}, pp.\,137-141)  of that precursor result, which  proceeds by extracting a   {\it martingale difference} subsequence. This crucial idea, which  establishes a powerful link to martingale theory and simplifies the arguments,    appears in this context for the first time  in \cite{K} (for related results, see   \cite{K1}).

\begin{theorem}
[\textsc{R\'ev\'esz} (1965)]
\label{Rev}
If the sequence $ \{f_n \}_{n \in \mathbb{N}}$ satisfies $\,\sup_{n \in \mathbb{N}} \mathbb{E} (f_n^2) < \infty\,,$   there exist a function $g \in \mathbb{L}^2$ and a subsequence $  \{f_{n_k}  \}_{k \in \mathbb{N}}\,,$  such that   
$ 
\, \sum_{k \in \mathbb{N}} a_k \big(f_{n_k} -g \big)\,
$ 
converges $\,\mathbb{P}-$a.e.,  for any sequence  $ \{a_k  \}_{k \in \mathbb{N}}  \subset \R$     with $\,\sum_{k \in \mathbb{N}} a^2_k < \infty$.
\end{theorem}

 It is clear   that this property of the subsequence $  \{f_{n_k}  \}_{k \in \mathbb{N}} $ is inherited by all {\it its} subsequences (just ``stretch out" the $a_k$'s accordingly, and fill out the gaps with zeroes).

In a   related development, \textsc{Delbaen \&   Schachermayer} (\cite{DS1},\,Lemma A1.1;\,\cite{DS2})   showed with very simple arguments that,  from every sequence $\{f_n \}_{n \in \mathbb{N}}$ of nonnegative, measurable functions,    a sequence of   convex combinations $\,g_n \in \text{conv}(f_n, f_{n+1}, \cdots ), ~ n \in \N\,$ of its elements can be extracted,  which converges $\mathbb{P}-$a.e.\,to a   measurable   $f_* : \Omega \to [0, \infty]$.    This result  was called ``a somewhat vulgar version of \textsc{Koml\'os}'s 
theorem" in \cite{DS2}, and  is implied 
  by  Theorem \ref{theorem3} below. Indeed,  
  convergence    for \textsc{Ces\`aro}  averages is much more precise than 
  for unspecified forward convex combinations.

In several contexts,    including optimization  treated via convex duality, nonnegativity is often no restriction at all, but rather  the natural setting      (e.g.,\,\cite{KS};\,\cite{LZ}; \cite{KS2};\,\cite{KK}, Chapter 3 and Appendix). Then, in the presence of convexity,   Lemma A1.1 in \cite{DS1}, or Theorem \ref{theorem3} here,    are very useful analogues of Theorem \ref{Kom}:  they lead  to limit functions $f_*$ in   convex sets  (such as the positive orthant in $   \mathbb{L}^0$, or the unit ball in $   \mathbb{L}^1$)   which are not compact in the usual sense, but   {\it  are}  ``convexly  compact"  as in \textsc{\v Zitkovi\'c} \cite{Z}.

\section{Result}
\label{sec3}

The purpose   of this note is to   prove with new and elementary tools the following version of Theorem \ref{Kom}, due to  \textsc{von\,Weizs\"acker}   \cite{vW}    
and   studied further     in \cite{Ta},   \S\,5.2.3 of  \cite{KS2}. 

\begin{theorem}
\label{theorem3}
Given a sequence $ \{f_n \}_{n \in \mathbb{N}}$ of  {\rm nonnegative}, measurable functions on a probability space $(\Omega, \mathcal{F}, \mathbb{P})$, there exist a measurable function $f_*:\Omega \to [0, \infty]$ and a subsequence $\big\{f_{n_k}\big\}_{k \in \mathbb{N}}$ of the original sequence, such that \eqref{02} holds.
\end{theorem}

Our     proof 
appears in Section \ref{sec5}; it is, we believe, not without methodological/pedagogical merit.  We observe that the result imposes no restriction whatsoever on the functions $f_1, f_2, \cdots$, apart from   measurability and nonnegativity. This comes at a price: the   function $f_*\,$,    constructed here carefully   in \eqref{5}--\eqref{8} below, can       take the value $+\infty$  on a set of positive measure. 




\section{Preparation}
\label{sec4}

We place ourselves in the setting of Theorem \ref{theorem3}. The arguments that follow often  necessitate    passing  to subsequences, and  to   diagonal subsequences, of a given   $ \{f_n \}_{n \in \mathbb{N}}$. To simplify typography, we denote  frequently such subsequences by the same symbols, $ \{f_n \}_{n \in \mathbb{N}}$. 

For each integer $k \in \mathbb{N}$, we introduce now the  truncated  functions
\begin{equation}
\label{3}
f_n^{(k)}\, :=\,  f_n \cdot  \mathbf{ 1}_{ \{  k-1 \le f_n <  k \}   }\,, \qquad n \in \mathbb{N}
\end{equation}
and note   the   partition of unity    
$\,
\sum_{k \in \mathbb{N}} f^{(k)}_n = f_n\,, ~ \forall ~ n \in \mathbb{N}\,.
$

\begin{lemma}
\label{lemma04}
For the sequence of functions $ \{f_n \}_{n \in \mathbb{N}}$ in Theorem \ref{theorem3}, there exists a subsequence, denoted by the same symbols  and such that, for every    $k\in \mathbb{N}$,     the functions of \eqref{3}  converge  to an appropriate measurable function $f^{(k)}:\Omega \to [0,\infty)\,,$ in the sense 
\begin{equation}
\label{4}
f_n^{(k)}  \, \xrightarrow[n \to \infty]{hC} \, f^{(k)}, \qquad \mathbb{P}-\text{a.e.}  
\end{equation}
For each fixed $k \in \N,$ this convergence holds also in $\mathbb{L}^1.$   
 \end{lemma}
 
 \noindent
 {\it Proof} (after \cite{Ch}, pp.\,145--146):   For arbitrary,   fixed     $k\in \mathbb{N}\,,$  the sequence $ \big\{f_n^{(k)} \big\}_{n \in \mathbb{N}}$ of \eqref{3} is bounded     in $\mathbb{L}^\infty,$ thus also in $\mathbb{L}^2$.   Theorem \ref{Rev} provides a function $ f^{(k)} \in \mathbb{L}^2$ and a   subsequence $  \{ f_{n_j}^{(k)}  \}_{j \in \mathbb{N}}$ of $  \{f_n^{(k)}   \}_{n \in \mathbb{N}}\,$, such that   $\,\sum_{j \in \mathbb{N}}    (f_{n_j}^{(k)} -f^{(k)}  )/ j\,$ converges $\,\mathbb{P}-$a.e.;   and as mentioned right after Theorem \ref{Rev}, this  is inherited by all subsequences of $  \{ f_{n_j}^{(k)}  \}_{j \in \mathbb{N}}$, and the \textsc{Kronecker}  Lemma (\cite{Du}, p.\,81)  gives   
 $$
 0= \lim_{J\to \infty} \frac{1}{J} \sum_{j=1}^J  \big( f_{n_j}^{(k)} - f^{(k)} \big)=\lim_{ J \to \infty} \frac{1}{J} \sum_{j=1}^J f_{n_j}^{(k)} - f^{(k)} ,\qquad \mathbb{P} -\hbox{a.e.}
 $$
We pass now  to a diagonal subsequence,  denoted  $ \big\{f_n  \big\}_{n \in \mathbb{N}}\,$ again,   and such that  \eqref{4} holds   for {\it every} $k \in \N\,.$ The last claim 
follows by the dominated convergence theorem. \qed

 \bigskip
With these ingredients, we introduce the measurable function $f:\Omega \to [0,\infty]$ via
\begin{equation}
\label{5}
f \,:=\, \sum_{k \in \mathbb{N}} f^{(k)}, \qquad \text{and consider the set} \quad A_\infty \,:= \,\{f=\infty\}.
\end{equation}
With the help of \textsc{Fatou}'s Lemma, and  the notation of (\ref{3})--(\ref{5}), Lemma \ref{lemma04} gives then
\begin{equation}
\label{6}
\varliminf_{N\to \infty} \frac{1}{N} \sum^{N}_{n=1} f_n \geq f\,, \qquad \mathbb{P} -\text{a.e.}
\end{equation}
\begin{equation}
\label{7}
\lim_{N\to \infty} \frac{1}{N} \sum^{N}_{n=1} f_n = \infty =f\,, \qquad \mathbb{P} -\text{a.e.} \quad \text{on} \quad A_\infty
\end{equation}
   for a suitable subsequence (denoted by the same symbols) of the original sequence $ \{f_n \}_{n \in \mathbb{N}}\,,$  and for all further subsequences of this subsequence.

The inequality in \eqref{6} can easily be strict. Consider, for instance, $f_n \equiv n\,,$ so that      $f_n^{(k)} =0$ holds in \eqref{3} for every fixed $k \in \mathbb{N}$ and all $n \in \mathbb{N}$ sufficiently large. We obtain  $f^{(k)} =0$ in \eqref{4}, thus $f=0$ in \eqref{5}; and yet $\frac{1}{N} \sum^N_{n=1} f_n \to \infty$ as $N\to \infty$.

\smallskip
 
This preparation allows us to formulate a   more technical and precise version of   Theorem \ref{theorem3},     Proposition \ref{prop05}   below, which   implies it. The   convention $ \,\infty \cdot 0=0$ is employed here, and  throughout.

\begin{proposition}
\label{prop05}
Fix a sequence $\{ f_n \}_{n \in \N}$ of nonnegative, measurable functions on the  probability space $(\Omega, \mathcal{F}, \mathbb{P})$, and recall the notation   \eqref{3}--\eqref{5}.  There exist then a subsequence, denoted again $\{f_n\}_{n\in \mathbb{N}}\,$,  and a  set $A \supseteq  A_\infty\,,$  
such that
\begin{equation}
\label{8}
f_n \, \xrightarrow[n \to \infty]{hC} \,f_* := \max \big(f, \, \infty \cdot \mathbf{1}_{A } \big)\,, \qquad \mathbb{P} -\text{a.e.}
\end{equation}
We have $A =A_\infty \,, $  thus also $f_* \equiv f ,$  when $\,\lim_{K\to \infty} \varlimsup_{n \to \infty} \, \mathbb{P} \big(f_n \ge K , f < \infty\big) =0\,$.  
\end{proposition}
  
This  last condition holds if    $\, \big\{ f_n \,\mathbf{1}_{ \{ f < \infty\}  } \big\}_{n \in \N}\,$ is bounded in $\mathbb{L}^0$,    i.e.,  $\,\lim_{K\to \infty} \sup_{n \in \N}  \, \P (f_n \ge K, f < \infty ) = 0\, .$    
A bit more stringently,        if not only $\{ f_n \}_{n \in \N}$ but also its solid,  convex hull in $\mathbb{L}^0_+,$ is bounded in $\mathbb{L}^0, $ then   $\{ f_n \}_{n \in \N}$ is bounded  in $\mathbb{L}
^1(\mathbb{Q})$ under some probability measure $ \mathbb{Q} \sim \mathbb{P},$ and thus $\mathbb{P} ( f < \infty) =1$ (e.g., Proposition A.11 in \cite{KK}). 
 Whereas, if $\{ f_n \}_{n \in \N}$ is  bounded   in $\mathbb{L}^1(\mathbb{P})$, i.e.,    $  \kappa := \sup_{n \in  \mathbb{N}} \mathbb{E} (f_n)< \infty  , $ then     $f$ in \eqref{5} is     integrable, since       $\,\mathbb{E} (f ) \le \kappa $ holds from \eqref{6} and \textsc{Fatou}.

\section{Proofs}
\label{sec5}

We shall need a couple of auxiliary results. First, and always with the notation of (\ref{3})--(\ref{5}), we note   the following consequence of monotone and dominated convergence.

\begin{lemma}
\label{lemma_2}
Suppose a set $\,D\subseteq \Omega \backslash A_\infty = \{f < \infty \}\,$ satisfies $\,\mathbb{E}\big(f \, \mathbf{1}_D\big) < \infty\,.$ 
Then, for any given  $\varepsilon \in (0,1),$  there exist   $K \in \N $   and a subsequence of the given sequence $\{ f_n \}_{n \in \N}$ such that  for it, and for any of its subsequences $($denoted again $\{ f_n \}_{n \in \N}),$  we have for arbitrary integers $L >K:$
\begin{equation}
\label{10}
\lim_{n \to \infty} \mathbb{E} \big[\,f_n \, \mathbf{ 1}_{   \{ K \le f_n <L \}\cap D }\, \big]=:\lim_{n \to \infty} \mathbb{E} \Big[\,f_n^{\,[K,L)}\, \mathbf{ 1}_D\, \Big] < \varepsilon\,. 
\end{equation}
\end{lemma}

 \smallskip
We are using  throughout  the notation   
\begin{equation}
\label{11}
f_n^{\,[K,L)} \,:=  \sum^L_{k=K+1} f_n^{(k)} = f_n \, \mathbf{1}_{[K,L)} (f_n)\,, \quad ~~~f_n^{\,[K,\infty)} \,:=\, \sum_{k \geq K+1} f_n^{(k)} = f_n  \, \mathbf{1}_{[K, \infty)} (f_n)   \,;
\end{equation}
  in an analogous manner  
$\,
f^{\,[K,L)}  :=   \sum^L_{k=K+1} f^{(k)}\,, ~   \,f^{\,[K,\infty)}  :=  \sum_{k \geq K+1} f^{(k)} 
\,,$  and   Lemma \ref{lemma04} gives 
\begin{equation}
\label{11a}
   f_n^{\,[K,L)}  \, \xrightarrow[n \to \infty]{hC} \,    f^{\,[K,L)}  
   \,,    \qquad   \hbox{both $\,\mathbb{P}-$a.e. and in $\,\mathbb{L}^1$.}  
\end{equation}

Secondly, we recall \eqref{7} and  observe the following dichotomy.   

\begin{lemma}
\label{lemma_3}
In the setting of   Proposition \ref{prop05}, consider any  measurable set $B \supseteq  \{f=\infty\} $  such that     the property  $\,f_n  \xrightarrow[n \to \infty]{hC} \infty\,$ of \eqref{7} holds $\,\mathbb{P}-$a.e.\,on $B$. 
Then, either

\medskip
\noindent
(i) there exist a set $\,C \supseteq B$ with $\mathbb{P} (C) > \mathbb{P} (B)$ and a subsequence, still denoted   $ \{f_n \}_{n \in \mathbb{N}}\,,$    with 
\begin{equation}
\label{9}
f_n  \xrightarrow[n \to \infty]{hC} \infty \qquad \text{valid} \quad \mathbb{P} -\text{a.e.} ~ \text{on} ~ \,C\,; \qquad \text{or,}
\end{equation}
(ii) the \textsc{Ces\`aro}  convergence   
$~f_n \xrightarrow[n \to \infty]{hC} f < \infty  ~~ \text{holds}~~  \mathbb{P} -a.e. ~ \text{on} ~\, \Omega \setminus B  \subseteq \{f< \infty\}  \,.$ 
\end{lemma}

Under {\it Case\,(ii),}    the set $B\supseteq A_\infty = \{f = \infty \}$ is  maximal  for the $\,\mathbb{P} -$a.e.\,\,property  $f_n \stackrel{hC}{\longrightarrow} \infty\,$: it cannot be ``inflated" to a set $C\supseteq B,$   which satisfies  \eqref{9} and  has bigger measure.  
This leads eventually to  Proposition \ref{prop05}, and thence to Theorem \ref{theorem3}.  

\smallskip
Before proving   these two results, we dispense with  the proof of Theorem \ref{Rev}; this is completely self-contained, and has nothing to do with Lemma \ref{lemma_2} or Lemma \ref{lemma_3}.

\subsection{Proof of Theorem \ref{Rev}}
\label{sec5f}

Because $\{f_n\}_{n \in \N}$ is bounded in $\mathbb{L}^2$, we can extract a subsequence that converges to some $g \in \mathbb{L}^2$ weakly in $\mathbb{L}^2$. Thus, it suffices    to prove the result for a sequence $\{g_n\}_{n \in \N}$ bounded in $\mathbb{L}^2$, and with $g_n \to 0$ weakly in $\mathbb{L}^2$. We take such a sequence, then, and approximate each $g_n$ by a {\it simple} function $h_n \in \mathbb{L}^2$ with  $\|g_n -h_n \|_2 \le 2^{-n}, ~\forall \,n \in \N.$ This gives, in particular,  
\begin{equation}
\label{R1}
\sum_{n \in \N} \, \big| g_n - h_n \big| < \infty\,, \quad \mathbb{P} -\hbox{a.e.\,;}\qquad h_n \to 0 \quad \hbox{weakly in } \mathbb{L}^2.
\end{equation}
We construct now,  by induction,  a sequence $1=n_1 < n_2 < \cdots\,$ of integers, such that
\begin{equation}
\label{R2}
 \big| \vartheta_k \big| < 2^{-k} \quad \text{holds }  ~ \mathbb{P} -\text{a.e., for } ~  ~\vartheta_k:= \E \big( h_{n_k} \, \big| \, h_{n_1} , \cdots, h_{n_{k-1}} \big)\,,~~k=2,3, \cdots ,
\end{equation}
as follows: The function $h_{n_1} = h_1$ is simple, thus so is  $ \,\E ( h_n | h_1) = \sum_{j=1}^J \gamma^{(n)}_j \mathbf{ 1}_{A_j} \,$ with $A_1, \cdots, A_J$ a partition of the space, and $\P (A_j)>0$,   $\gamma^{(n)}_j := \big(1 / \P (A_j) \big) \cdot \E \big( h_n \,\mathbf{ 1}_{A_j} \big)\,.$ This  last  expectation tends to zero as $n \to \infty$  from \eqref{R1}, for every fixed $j$; so we can choose $n_2 > n_1 =1$ with $ \big|\gamma^{(n_2)}_j \big| < 2^{-2},$ for
 \noindent
 $ j=1, \cdots, J$; i.e., $\big|\vartheta_2 \big|< 2^{-2},$ $\P-$a.e. Clearly, we can keep repeating  this argument  since, at each stage, $\big( h_{n_1}  ,   \cdots, h_{n_{k-1}} \big)$  generates a finite partition of the space; and this way we  arrive at \eqref{R2}. 

 \newpage
 
The sequence $\{ h_n \}_{n \in \N}$ is bounded in $\mathbb{L}^2$,  thus   so   is  the martingale
$ 
X_n  := \sum_{k=0}^n a_k \big(h_{n_k} - \vartheta_k \big)\,, ~ n \in \N_0\,,
$ 
for any       $ \{a_n  \}_{n \in \mathbb{N}_0} \subset \R  $    with $\sum_{n \in \mathbb{N}} a^2_n < \infty$. Martingale convergence theory (\cite{Du}, p.\,236) shows that   the series $\sum_{k \in \N}  a_k \big(h_{n_k} - \vartheta_k \big)$ converges $\P-$a.e. But we have also $\sum_{k \in \N}    \big(  \big| \vartheta_k \big| + \big| g_{n_k} - h_{n_k} \big| \big)< \infty, $ $\,\P-$a.e.\,\,from \eqref{R1}--\eqref{R2}, and   deduce that   
$ 
  \sum_{k \in \mathbb{N}} a_k \,g_{n_k}   
$ 
converges $\,\mathbb{P}-$a.e., the claim of the theorem.  \qed

\subsection{Proof of Lemma \ref{lemma_2}}
\label{sec5a}

Let us call {\it ``Lemma \ref{lemma_2}$^{\,\dagger}$"} the same    statement as that of Lemma \ref{lemma_2}, except that \eqref{10} is now replaced by 
\begin{equation}
\label{10a}
\forall~ L=K+1, K+2, \cdots \,:~~~  \mathbb{E} \Big[\,f_n^{\,[K,L)}\, \mathbf{ 1}_D\, \Big] < \varepsilon\,, ~\text{for all but finitely many}~ n \in \N .
\end{equation}
{\it Claim: Lemma \ref{lemma_2}$^{\,\dagger}$ implies Lemma \ref{lemma_2}.} Let a subsequence of the original $\{ f_n \}_{n \in \N}$ be given (denoted   $\{ f_n \}_{n \in \N}$ again), along with   arbitrary $\varepsilon \in (0,1)$. Lemma \ref{lemma_2}$^{\,\dagger}$ guarantees the existence of $K \in \N$, depending on $\varepsilon$ and the subsequence, such that \eqref{10a} holds for all integers $L \ge K+1$.

Choose $L=K+1$ first. From Lemma \ref{lemma_2}$^{\,\dagger}$ and \textsc{Bolzano-Weierstrass},   (the current) $\{ f_n \}_{n \in \N}$ has a  subsequence for which the expectation in \eqref{10a} converges, with limit $\le \varepsilon / 2.$ Now choose $L=K+2$ and a  subsequence of the last subsequence, for which the expectation in \eqref{10a}   converges and has  limit $\le \varepsilon / 2.$ Continuing in this manner, then   diagonalizing, we obtain a subsequence that satisfies \eqref{10a}. 

\bigskip
\noindent
{\it Proof of Lemma \ref{lemma_2}$^{\,\dagger}$.} We   argue by contradiction, assuming that $\{ f_n \}_{n \in \N}$ has a subsequence for which Lemma \ref{lemma_2}$^{\,\dagger}$ fails. Then there     exists an $\varepsilon \in (0,1)$ with the property that, for every subsequence of $\{ f_n \}_{n \in \N}$  and every $K \in \mathbb{N}$, there exists an integer $L>K$ such that  
\begin{equation}
\label{11b}
\mathbb{E} \bigg[\sum^L_{k=K+1}  f^{(k)}_n \, \mathbf{1}_D \bigg] = \mathbb{E} \Big(f_n^{[K,L)} \,  \mathbf{1}_D \Big) \geq \varepsilon
\end{equation}
holds for infinitely many integers $n \in \mathbb{N}$. But this means that there is a subsequence, again denoted by $\{f_n\}_{n \in \mathbb{N}}\,,$   {\it along which we have \eqref{11b}   for every $n \in \mathbb{N};$} and, as a result,   also
\begin{equation}
\label{11c}
\mathbb{E} \bigg[\sum^L_{k=K+1} \Big(\frac{1}{N} \sum^N_{n=1} f_n^{(k)}\Big) \, \mathbf{1}_D \bigg] \geq \varepsilon\,, \qquad \forall ~ n \in \N\,.
\end{equation}
Now  all the truncated functions $f_n^{(k)}$ as in \eqref{3}, for $k=K+1, \dots, L$ and $n \in \mathbb{N}$, take values on the ``Procrustean bed" $\{ 0 \} \cup [K,L)$; and    $\, \lim_{N\to\infty} \frac{1}{N}  \sum^N_{n=1} f_n^{(k)} = f^{(k)}\,$ holds $\mathbb{P} -$a.e., 
for the selected subsequence and all its subsequences,  on account of Lemma \ref{lemma04}. Thus, $\,
\mathbb{E} \big[\sum^L_{k=K+1} f^{(k)} \, \mathbf{1}_D \big] \geq \varepsilon
\,$ from bounded convergence and \eqref{11c}; and    the nonnegativity of these $f^{(k)}$'s implies   also
\begin{equation}\label{11d}
\mathbb{E} \bigg(\sum_{k\geq K+1} f^{(k)} \, \mathbf{1}_D \bigg) = \mathbb{E} \Big(f^{\,[K, \infty)} \, \mathbf{1}_D \Big) \geq \varepsilon\,, \qquad \forall ~~K \in \mathbb{N}\,.
\end{equation}

The  nonnegativity gives also $\,\lim_{K\to \infty} \uparrow \sum^K_{k=1} f^{(k)}\, \mathbf{1}_D =f\, \mathbf{1}_D\,$, both $\mathbb{P} -$a.e.\,and in $\mathbb{L}^1$. Since $ \mathbb{E}\big(f \, \mathbf{1}_D\big) < \infty $  by assumption,   $ \mathbb{E} \big[f^{\,[K, \infty)} \, \mathbf{1}_D\big] <  \varepsilon/2\,$ holds for   all  $K \in \mathbb{N}$   large enough. But this contradicts \eqref{11d}, and we are done.  \qed

\subsection{Proof of Lemma \ref{lemma_3}}
\label{sec5b}

 \smallskip

We start by fixing   $j \in \mathbb{N}$  and distinguishing two contingencies,  with the definitions  
\begin{equation}
\label{12}
D_j := \{f \leq j\} \backslash B
\,, \qquad 
E_n^{[K, \infty)} \,:=\, \big\{f_n^{\,[K,\infty)} \geq K \big\} \cap  D_j   \,=\, \big\{f_n  \geq K \big\} \cap  D_j    \,,
\end{equation}
\begin{equation}
\label{12too}
\alpha \,:=\,  \lim_{K\to \infty} \varlimsup_{n \to \infty} \mathbb{P} \big(E_n^{[K, \infty)}\big) \,:
\end{equation}
\newpage

\noindent
  {\it Contingency ~I:} $~\alpha > 0\,.$

\smallskip
\noindent
 {\it Contingency II:} $~\alpha = 0\,.$

\medskip
\noindent
$\bullet~$ 
Under {\bf Contingency\,I}\,, we pass to a subsequence $\{f_n\}_{n \in \mathbb{N}}$ with $\mathbb{P} \big(E_n^{\,[n^2, \infty)}\big) \geq  \alpha /2 \,,$   $\forall ~n \in \mathbb{N}\,$; and consider   indicators $\,g_n := \mathbf{1}_{ E_n^{\,[n^2, \infty)}},  ~n \in \mathbb{N}\,,$ all of them supported on the set $\Omega \setminus B$. Arguing as in Lemma \ref{lemma04} we obtain 
 a subsequence, still denoted    $\{g_n\} _{n \in \mathbb{N}}$,  with 
$ \,
g_n \xrightarrow[n \to \infty]{hC} g\,,$ $  \mathbb{P} -\hbox{a.e.,}
$ 
for some   $g:\Omega \to [0,1]\,$ with $\{g >0\} \subseteq \Omega \backslash   B \,$   and $\mathbb{E}(g) \geq \alpha /2$  by   bounded convergence.   

 \smallskip
Thus,   $f_n \xrightarrow[n \to \infty]{hC} \infty\, $  holds  $\,    \mathbb{P} -$a.e.  on   $\,\{g >0\}\,$. This set   
has  $\,\mathbb{P} \big( g>0\big) = \mathbb{E} [\mathbf{1}_{\{g>0\}} ] \geq \mathbb{E}(g) \geq \alpha/2 \,;$    we are   under {\it Case  (i)} of Lemma \ref{lemma_3}, with $C := \{g >0\} \cup B$ and $\mathbb{P}(C) > \mathbb{P}(B)$.

\medskip
\noindent
$\bullet~$
Now we pass to  {\bf Contingency\,II}\,.     We fix 
$\varepsilon >0,$   $D_j=\{f \leq j\}  \backslash B
$, and apply Lemma \ref{lemma_2} with this $D_j$ to construct inductively a subsequence $\big\{ n_m \big\}_{m \in \N}\,,$ along with sequences  $\big\{ K_m \big\}_{m \in \N}\,, $ $\big\{ L_m \big\}_{m \in \N}\, $ of integers increasing to infinity and such that 
\begin{equation}
\label{A.8}
  \mathbb{P} \big(E_{n_m}^{\,[L_{m }, \infty ) }\big) = \mathbb{P} \big(\big\{ f_{n_m} \ge L_m \big\} \cap D_j \big) < 2^{- m } 
\end{equation}
\begin{equation}
\label{A.9}
  \mathbb{E} \Big[\, f_{n_p}^{\,[K_m, L_p)} \, \mathbf{ 1}_{D_j} \, \Big] 
  < 2^{- m }\,,\qquad \forall ~~p = m, m+1, \cdots 
\end{equation}
hold for every $ m \in \N$. With the  choice (\ref{A.8}), the sequences $\, \big\{ f_{n_m} \cdot \mathbf{ 1}_{D_j} \big\}_{m \in \N}\,$ and $\, \big\{ f_{n_m}^{\,[0, L_m)} \cdot \mathbf{ 1}_{D_j} \big\}_{m \in \N}\,$ are equivalent in the sense introduced in section \ref{sec2}, as the probability of their respective general terms being different is bounded from above by $2^{- m }$. We claim that 
\begin{equation}
\label{A.14}
f_{n_m} \cdot \mathbf{ 1}_{D_j}\,  \xrightarrow[m \to \infty]{hC} \,f  \cdot \mathbf{ 1}_{D_j}\,,\quad   \mathbb{P} -\hbox{a.e.;}
\end{equation}
and in view of the previous statement, this amounts to 
\begin{equation}
\label{A.10}
f_{n_m}^{\,[0, L_m)} \cdot \mathbf{ 1}_{D_j}\,  \xrightarrow[m \to \infty]{hC} \,f  \cdot \mathbf{ 1}_{D_j}\,,\quad   \mathbb{P} -\hbox{a.e.}
\end{equation}

To prove (\ref{A.10}), we start by observing that the sequence $\, \big\{ f_{n_m}^{\,[0, L_m)} \cdot \mathbf{ 1}_{D_j} \big\}_{m \in \N}\,$ is {\it uniformly integrable}, thus bounded in $\mathbb{L}^1$, as 
$$
\sup_{p \in \N \atop p \ge m}  \,   \mathbb{E} \Big[\, f_{n_p}^{\,[0, L_p)} \, \mathbf{ 1}_{D_j} \cdot \mathbf{ 1}_{ \big\{ f_{n_p}^{\,[0, L_p)} \ge K_m \big\} }\, \Big] 
  < 2^{- m }
$$
holds on account of (\ref{A.9}) for every $m \in \N\,$.  Theorem \ref{Kom} gives   an integrable function $h : \Omega \to [0, \infty)\,$ with 
\begin{equation}
\label{A.11}
f_{n_m}^{\,[0, L_m)} \cdot \mathbf{ 1}_{D_j}\,  \xrightarrow[m \to \infty]{hC} \,h  \cdot \mathbf{ 1}_{D_j}\,,\quad   \mathbb{P} -\hbox{a.e.,}
\end{equation}
and we need to argue that this $h$ agrees with $f$ from (\ref{5}), $\mathbb{P}-$a.e.\,on $D_j$.

Indeed, for every $K \in \N$ and all $m$ large enough,   
$\,
\sum_{k=1}^K \, f_{n_m} \, \mathbf{ 1}_{  \{ k-1 \le f_{n_m} < k  \} } \,=\, f_{n_m}^{\,[0,K)} \,\le \, f_{n_m}^{\,[0,L_m)}\, 
 $
holds, therefore $\, \sum_{k=1}^K f^{(k)} \cdot \mathbf{ 1}_{D_j} \le h \cdot \mathbf{ 1}_{D_j}\,$ by letting $m \to \infty$, on account of (\ref{A.11}) and Lemma \ref{lemma04}. Passing now to the limit as $K \to \infty$ and recalling (4.3), we arrive at 
\begin{equation}
\label{A.13}
f  \cdot \mathbf{ 1}_{D_j}\, \le \, h  \cdot \mathbf{ 1}_{D_j}\,,\quad   \mathbb{P} -\hbox{a.e.}
\end{equation}

To obtain the inequality in the reverse direction, we take expectations. From   (\ref{A.11}) and uniform integrability, we have 
$\,
 \mathbb{E} \Big[\, f_{n_m}^{\,[0, L_m)} \cdot \mathbf{ 1}_{D_j} \, \Big] \,  \xrightarrow[m \to \infty]{hC}\,  \mathbb{E} \big[\, h  \cdot \mathbf{ 1}_{D_j} \, \big]\,,$ therefore also 
 $$
  \mathbb{E} \big[\, h  \cdot \mathbf{ 1}_{D_j} \, \big]\,=  \lim_{M  \to \infty \atop K  \to \infty} \frac{1}{M}\, \mathbb{E} \bigg[\, \sum_{m=1}^M f_{n_m}^{\,[0, L_m \wedge K)} \cdot \mathbf{ 1}_{D_j} \, \bigg] \,=  \lim_{M  \to \infty \atop K  \to \infty} \frac{1}{M}\, \mathbb{E} \bigg[\, \sum_{m=1}^M \sum_{k=1}^{L_m \wedge K} f_{n_m}^{(k)} \cdot \mathbf{ 1}_{D_j} \, \bigg]~~~~~~~~~~~~
 $$
 $$
~~~~~~~~~~~~~~ \le \lim_{M  \to \infty \atop K  \to \infty} \frac{1}{M}\, \mathbb{E} \bigg[\, \sum_{m=1}^M \bigg( \sum_{k=1}^{  K} f_{n_m}^{(k)} \bigg) \cdot \mathbf{ 1}_{D_j} \, \bigg]\,=\,\lim_{K \to \infty } \, \mathbb{E} \bigg[\, \sum_{k=1}^{K } \, f^{(k)}   \cdot \mathbf{ 1}_{D_j} \, \bigg]\,\le\, \mathbb{E} \big[\, f  \cdot \mathbf{ 1}_{D_j} \, \big]
 $$
\newpage
\noindent
from  Lemma \ref{lemma04}. In conjunction with (\ref{A.13}), this shows $\,f  \cdot \mathbf{ 1}_{D_j}\, =\, h  \cdot \mathbf{ 1}_{D_j}\,,~~   \mathbb{P} -\hbox{a.e.},$ as claimed; and on account of (\ref{A.11}) it establishes (\ref{A.10}), thus (\ref{A.14}) as well.

\smallskip
The final step is to let $ \,j\to \infty\,$: we do this again by  extracting subsequences, successively for each $j \in \N\,$, then  passing   to a diagonal subsequence. We    obtain then  \eqref{A.14} with $D_j$ replaced   by the set  $D:= \bigcup_{j \in \N}D_j =\{ f < \infty \}  \backslash B
,$   and deduce that we are in  {\it Case\,(ii)} of Lemma \ref{lemma_3}. \qed

\subsection{Proofs of Proposition \ref{prop05} and Theorem \ref{theorem3}}
\label{sec5c}

On the strength of Lemma \ref{lemma_3} we construct,  by exhaustion or transfinite induction arguments  and as long as we are under   the dispensation of its {\it Case\,(i),} an increasing sequence $  B\subseteq  B_1 \subseteq B_2 \subseteq  \dots$ of sets as postulated there, whose union $B_\infty := \bigcup_{j \in \mathbb{N}} B_j \supseteq   B \supseteq \{f=\infty\}$ is maximal with the property \eqref{9} for an appropriate subsequence. But   maximality means that, on the complement  $\Omega \backslash B_\infty$ of this set,  we must be in  the realm   of {\it Case\,(ii)} in Lemma \ref{lemma_3}. This   establishes   the first claim of    Proposition \ref{prop05} with $A = B_\infty \supseteq  \{f=\infty\} \,$,    thus also   Theorem \ref{theorem3}.     

 For the second claim of the Proposition,  we note that   equality holds right above, that is, $  B_\infty =  \{f=\infty\}  ,$     if we  are under  Contingency II (i.e.,    $\alpha =0$) in \S\,\ref{sec5b} (proof of Lemma \ref{lemma_3}) and with $B   = \{f=\infty\}\,$ in (\ref{12}); a sufficient condition for this, is $\, \lim_{K\to \infty} \varlimsup_{n \to \infty} \mathbb{P} (f_n \ge K, f < \infty ) =0$.           The claim now follows.       \qed


   \end{document}